\documentclass [11pt]{article}
\usepackage{latexsym,amsfonts}
\begin{document}
\bibliographystyle{plain}

\newtheorem{defi}{Definition}[section] 
\newtheorem{theo}{Theorem}[section]
\newtheorem{prop}{Proposition}[section] 
\newtheorem{lemm}{Lemma}[section]     
\newtheorem{conj}{Conjecture}[section] 
\newtheorem{note}{Note}[section] 
\newtheorem{coro}{Corollary}[section]    
\newtheorem{ques}{Question}[section] 
\newtheorem{rema}{Remark}[section]

\newcommand{\fas}{\ensuremath{f_{\rm {A, sol}}(n)}} 
\newcommand{\ms}{\ensuremath{M_{\rm{s,p',A}}(G)}} 
\newcommand{\fbh}{\ensuremath{f({\mathbb Z}_{p^\beta}H, \alpha)}}   
\newcommand{\fah}{\ensuremath{f({\mathbb Z}_{p^\alpha}H, \alpha)}}   
\newcommand{\zah}{\ensuremath{{\mathbb Z}_{p^\alpha}H}} 
\newcommand{\zbh}{\ensuremath{{\mathbb Z}_{p^\beta}H}} 
\newcommand{\zih}{\ensuremath{{\mathbb Z}_{p^\alpha}H_i}} 
\newcommand{\zpq}{{\mathbb Z}_pQ} 
\newcommand{\zph}{{\mathbb Z}_pH} 
\newcommand{\phx}{\ensuremath{{\phi}_{\scriptscriptstyle X}}}   
\newcommand{\jpq}
{\ensuremath{{\mathfrak A}_p{\mathfrak A}_q \vee {\mathfrak A}_q{\mathfrak A}_p}} 
\newcommand{\pqr}
{\ensuremath{{\mathfrak A}_p{\mathfrak A}_q{\mathfrak A}_r}}
\newcommand{\jpqr}
{\ensuremath{ \bigvee_{\pi}\left\{{\mathfrak A}_{\pi(p)}{\mathfrak A}_{\pi(q)}
{\mathfrak A}_{\pi(r)}\right\}}}
\newcommand{\fu}{\ensuremath{f_{\mathfrak U}(n)}}
\newcommand{\fv}{\ensuremath{f_{\mathfrak V}(n)}}
\newcommand{\fw}{\ensuremath{f_{\mathfrak W}(n)}} 
\newcommand{\fvp}{\ensuremath{f_{\mathfrak V}(p^{\alpha}q^{\beta}r^{\gamma})}} 
\newcommand{\fxp}{\ensuremath{f_{\cal X}(p^{{\alpha}_1}q^{{\beta}_1}r^{{\gamma}_1})}} 
\newcommand{\fy}{\ensuremath{f_{\cal Y}(p^{\alpha}_2q^{\beta}_2)}} 
\newcommand{\fbx}{\ensuremath{f_{\bar {\cal X}}(p^{\alpha}_1q^{\beta}_1)}} 
\newcommand{\fby}{\ensuremath{f_{\bar {\cal Y}}(p^{\alpha_2}q^{\beta_2})}} 

\newcommand{\lav}{\ensuremath{a = \limsup\left \{\frac{\log f_{\mathfrak V}(n)}{\mu (n) \log n}\right \} }}
\newcommand{\lau}{\ensuremath{a = \limsup\left \{\frac{\log f_{\mathfrak U}(n)}{\mu (n) \log n}\right \} 
}}

\newcommand{\law}{\ensuremath{a = \limsup\left \{\frac{\log f_{\mathfrak W}(n)}{\mu (n) \log n}\right \} 
}}

\title{On irreducibility of induced modules and an adaptation of the Wigner--Mackey method of little groups} 

\author{Geetha Venkataraman \\
St. Stephen's College \\
University of Delhi \\
Delhi-110007 \\
India \\
\vspace{.05in}
email: geevenkat@gmail.com}
\date{\today}            
\maketitle

\noindent {\bf Abstract}
This paper deals with sufficiency conditions for irreducibility of certain induced modules. We also construct irreducible representations for a group $G$ over a field  ${\mathbb K}$ where the group $G$ is a semidirect product of a normal abelian subgroup $N$ and a subgroup $H$.  The main results are proved with the assumption that ${\rm char} \, {\mathbb K}$ does not divide $|G|$ but there is no assumption made of ${\mathbb K}$ being algebraically closed.

\vspace{.15in}

\vspace{.15in}

\section{Introduction}  

In this paper $G$ is a finite group and $\mathbb K$ is an arbitrary field unless stated otherwise. We shall also assume that the ${\mathbb K}G$-modules referred to here are finite dimensional over ${\mathbb K}$ and that these modules are left ${\mathbb K}G$-modules. For proofs of well known theorems that are used here and terminology see \cite{cr62} or \cite{gk90}.

Mackey \cite{mgw51} proved results about necessary and sufficient conditions for the irreducibility of induced modules. A crucial assumption was that the field over which the representations occur be algebraically closed.

We show in this paper that when the condition of  the base field being an algebraically closed field is dropped then the sufficiency condition for the induced module to be irreducible still holds. We prove the following.

\vspace{.05in}
\noindent{\bf Theorem 3.1}
{\em Let $G$ be a finite group and let ${\mathbb K}$ be a field such that ${\rm char}\,{\mathbb K}$ does not divide $|G|$. Let $H$ be a subgroup of $G$ and let $L$ be an irreducible ${\mathbb K}H$-module. For $x \in G$, let $H^{(x)} = xHx^{-1} \cap H$. If for all $x \not\in H$, the ${\mathbb K}H^{(x)}$-modules, $L_{H^{(x)}}$ and $(x \otimes L)_{H^{(x)}}$ are disjoint, {\rm(}that is, they have no composition factors in common{\rm)} then $L^G$ is an irreducible ${\mathbb K}G$-module.}
\vspace{.05in}

The other topic dealt with in this paper follows the same approach as the method of {\it little groups} of Wigner and Mackey (see \cite{jps77} pp 62--63). Let $G$ be a finite group which is a semidirect product of a normal abelian subgroup $N$ by $H$. Let ${\mathbb K}$ be an algebraically closed field such that ${\rm char} \, {\mathbb K}$ does not divide $|G|$. Then the method of {\it little groups} shows how the irreducible representations of $G$ over ${\mathbb K}$ can be constructed from those of certain subgroups of $H$. 

Note that since ${\mathbb K}$ is algebraically closed and since $N$ is abelian, all its irreducible representations are of degree 1. Mackey's irreducibility criterion is used in the method of little groups to prove the irreducibility of certain induced modules and so ${\mathbb K}$ being algebraically closed is essential in Wigner--Mackey method of little groups for constructing the irreducible representations of $G$.

Using Theorem 3.1 stated above, we are able to follow the same proof as in the Wigner--Mackey method of little groups to get a classification of irreducible representations of $G$ except that the field ${\mathbb K}$ is no longer required to be algebraically closed. However we do need the condition that all irreducible representations of $N$ are of degree 1. If we drop this last condition then we can classify all irreducible ${\mathbb K}G$-modules which have a one-dimensional composition factor when restricted to $N$. These results are presented as Theorems 4.1 and 4.2.

This paper is organised as follows.
In Section 2, we present the preliminary results required for the proofs of our main results. Some well know theorems that are required for their proofs are also stated in this section. In the third section we present our main result on irreducibility of induced modules and the corollaries following from it. Section 4 deals with applications of Theorem 3.1, namely, the results following from our adaptation of the Wigner--Mackey method of little groups. 

\section{Intertwining number}

We begin by establishing some definitions and notations.

For an arbitrary associative ring $R$, and $R$-modules $M$ and $N$ we denote by ${\rm Hom}_R(M, N)$ the additive group consisting of all $R$-homomorphisms from $M$ to $N$. 

\begin{defi}
Let $M$ and $N$ be ${\mathbb K}G$-modules. Then ${\rm Hom}_{{\mathbb K}G}(M, N)$ is a vector space over ${\mathbb K}$, and its dimension is called the intertwining number $i(M, N)$ of $M$ and $N$.
\end{defi}

Some basic properties of the intertwining number are stated below. All parts are easy to prove.
\begin{rema}
Let $M$, $N$, $M_i$, $N_i$, $i = 1,2$ be ${\mathbb K}G$-modules. Then
\begin{enumerate}
\item[\rm{(a)}]$i(M_1 \oplus M_2, N) \, = \, i(M_1, N) + i(M_2, N)$.
\item[\rm{(b)}]$i(N, M_1 \oplus M_2)\,  = \, i(N, M_1) + i(N, M_2)$.
\item[\rm{(c)}]If $M$ and $N$ are completely reducible as ${\mathbb K}G$-modules, then $i(M, N) \, = \, i(N, M)$.
\item[\rm{(d)}]Let $M$ and $N$ be irreducible ${\mathbb K}G$-modules. Then $i(M, N) \, = \, 0$ if and only if $M \not \cong N$ as ${\mathbb K}G$-modules.
\item[\rm{(e)}]If $M \not = 0$ then $i(M, M) \not = 0$. Further if dimension of $M$ over ${\mathbb K}$ is 1, then $i(M, M) = 1$.
\item[\rm{(f)}]Let $M$ and $N$ be completely reducible ${\mathbb K}G$-modules. Then $M$ and $N$ are disjoint if and only if $i(M, N) \, = \, 0$. {\rm(}The modules $M$ and $N$ are said to be disjoint if they have no composition factor in common.{\rm)}
\end{enumerate}
\end{rema}

We use two well known theorems in the proof of our main results and the statements of these are presented below.
 
 \vspace{.05in}
\noindent
The following is the statement of the Frobenius Reciprocity Theorem (FRT) for ${\mathbb K}G$-modules. (For a proof of this result see \cite{crv181}, pp 232--233.) 

 \vspace{.05in}
\noindent{\bf Theorem: (FRT)}{  }{\em Let $H$ be a subgroup of $G$. Let $V$ be a ${\mathbb K}G$-module and let $W$ be a ${\mathbb K}H$-module. Then
$$
{\rm Hom}_{{\mathbb K}G}(W^G, V) \cong {\rm Hom}_{{\mathbb K}H}(W, V_H) \: .
$$}

\vspace{.05in}
\noindent
The statement of the Intertwining Number Theorem (INT) is presented below. (For a proof of this result see \cite{cr62}, p 327.) 

 \vspace{.05in}
\noindent{\bf Theorem: (INT)}{  }{\em Let $H_1$ and $H_2$ be subgroups of $G$ and let $L_i$ be ${\mathbb K}H_i$-modules for $i=1, 2$. Let $(x, y) \in G \times G$. Set $H^{(x,\, y)} = xH_1x^{-1} \cap yH_2y^{-1}$. Further let ${L_1}^{(x)} := x \otimes L_1 \subseteq {L_1}^G$ and let 
${L_2}^{(y)} := y \otimes L_2 \subseteq {L_2}^G$. Then 
\begin{enumerate}
\item[\rm{(i)}]${L_1}^{(x)}$ and ${L_2}^{(y)}$ are ${\mathbb K}H^{(x,\, y)}$-modules.  
\item[\rm{(ii)}]The intertwining number $i({L_1}^{(x)}, {L_2}^{(y)})$ of the ${\mathbb K}H^{(x,\, y)}$-modules depends only on the $(H_1, H_2)$-double coset $D$ to which $x^{-1}y$ belongs and will be denoted as $i(L_1, L_2, D)$.
\item[\rm{(iii)}]$i({L_1}^G, {L_2}^G) = \sum_{D}i(L_1, L_2, D)$ where the sum is taken over all $(H_1, H_2)$-double cosets $D$ in $G$.
\end{enumerate}
}

\section{Irreducibility of induced modules}

We prove a result required for the proof of our main theorem.

\begin{prop}
Let $H$ be a subgroup of $G$. Let ${\mathbb K}$ be a field such that ${\rm char}\,{\mathbb K}$ does not divide the order of $G$. Let $L$ be an irreducible ${\mathbb K}H$-module. If $i(L^G, L^G) =
i(L, L)$ then $L^G$ is an irreducible ${\mathbb K}G$-module.
\end{prop}

\noindent{\bf Proof}{   }Let $i(L^G, L^G) = i(L, L)$ and let $L^G  = \oplus_{j=1}^s M_j$ be the decomposition of the completely reducible ${\mathbb K}G$-module $L^G$ into irreducible ${\mathbb K}G$-modules. Then ${(L^G)}_H  = \oplus_{j=1}^s {(M_j)}_H$ as ${\mathbb K}H$-modules.

By FRT,  for $j=1, \ldots, s$, we have,
$$
{\rm Hom}_{{\mathbb K}G}(L^G, M_j) \cong {\rm Hom}_{{\mathbb K}H}(L, {(M_j)}_H) \: .
$$
So, $i(L^G, M_j)= i(L, {(M_j)}_H)$ and hence $i(L, {(M_j)}_H) \not = 0$ for $j=1, \ldots, s$.

Therefore using Remark 2.1 and the above fact, for any $j=1, \ldots, s$, there is an irreducible ${\mathbb K}H$-submodule of $M_j$ which is isomorphic to $L$. Hence for all $j=1, \ldots, s$,
\begin{eqnarray}
i(L, {(M_j)}_H) & \geq & i(L, L) \: .
\end{eqnarray}

By FRT we have,
$$
{\rm Hom}_{{\mathbb K}G}(L^G, L^G) \cong {\rm Hom}_{{\mathbb K}H}(L, {(L^G)}_H) \: .
$$
Thus by the above isomorphism, Remark 2.1, part (b) and (1), we get 
\begin{eqnarray*}
i(L, L) & = &i(L^G, L^G)\\
        & = &i(L, {(L^G)}_H)\\
& = & \sum_{j=1}^s i(L, {(M_j)}_H) \\
 & \geq & s \ i(L, L) \:.
\end{eqnarray*}
Consequently we must have $s=1$, and so $L^G$ is irreducible.\ $\Box$

Now we are in a position to present our main result.

\begin{theo}
Let $G$ be a finite group and let ${\mathbb K}$ be a field such that ${\rm char}\,{\mathbb K}$ does not divide $|G|$. Let $H$ be a subgroup of $G$ and let $L$ be an irreducible ${\mathbb K}H$-module. For $x \in G$, let $H^{(x)} = xHx^{-1} \cap H$. If for all $x \not\in H$, the ${\mathbb K}H^{(x)}$-modules, $L_{H^{(x)}}$ and $(x \otimes L)_{H^{(x)}}$ are disjoint then $L^G$ is an irreducible ${\mathbb K}G$-module.
\end{theo}

\noindent{\bf Proof}{   }By INT we have that
\begin{eqnarray}
i({L}^G, {L}^G) & = & \sum_{D}i(L, L, D)
\end{eqnarray} 
where the sum is taken over all $(H, H)$-double cosets $D$ in $G$. 

Note that $i(L, L, D) = i(L^{(x)}, L^{(y)})$ for some $x^{-1}y \in D$, where $L^{(x)}$ and $L^{(y)}$ are $H^{(x,\, y)}$ modules over ${\mathbb K}$ and $H^{(x,\, y)} =xHx^{-1} \cap yHy^{-1}$. If we take $y=1$, then $H^{(x,\, y)} = H^{(x)}$. So $L^{(y)} = L_{H^{(x)}}$, $x^{-1} \in D$ and $i(L, L, D) = i(L^{(x)}, L_{H^{(x)}})$. 

For the double coset $D=H$, we can take $x=1$ and in this case we get $H^{(x)} = H$ and so 
\begin{eqnarray}
i(L, L, H) & = & i(L, L) \: .
\end{eqnarray}

For any other double coset $D$ different from $H$, if $x^{-1} \in D$, then $x^{-1} \not \in H$ or equivalently $x \not \in H$. But then we are given that ${(L^{(x)})}_{H^{(x)}}$ and $L_{H^{(x)}}$ are disjoint as ${\mathbb K}H^{(x)}$-modules. So for $D \not = H$, using Remark 2.1, we can show easily that
\begin{eqnarray}
i(L, L, D) & = & i({(L^{(x)})}_{H^{(x)}}, L_{H^{(x)}}) \nonumber \\
& = & 0 \: .
\end{eqnarray}

Substituting the values of the intertwining numbers given in (3) and (4), in equation (2), we have

\begin{eqnarray*}
i(L^G, L^G) & = & i(L, L, H) + \sum_{D \not = H}i(L, L, D)\\
& = & i(L, L) \: .
\end{eqnarray*}

So by Proposition 3.1 the induced module, $L^G$ is irreducible as a ${\mathbb K}G$-module. \hspace{.75in} $\Box$

Before we present the corollaries arising from the above theorem, we have the following remark.
\begin{rema}

From the proof of the above theorem it is easy to see in fact that $i(L^G, L^G) = i(L, L)$ if and only if for all $x \not\in H$, the ${\mathbb K}H^{(x)}$-modules, $L_{H^{(x)}}$ and $(x \otimes L)_{H^{(x)}}$ are disjoint.

On the other hand, while the equivalent conditions given above are sufficient for $L^G$ to be irreducible, they are not necessary. A simple example to show that this condition is not necessary is the following. 

Let $G$ be the cyclic group of order 4 and let $H$ be the unique subgroup of order 2. Let ${\mathbb K}$ be the field of three elements and $L$ be the only non-trivial irreducible ${\mathbb K}H$-module of dimension 1. Then we can show that $L^G$ is irreducible and that $i(L^G, L^G) = 2 = 2 \, i(L,L)$.
\end{rema}

We have the following corollaries to the above theorem.
\begin{coro}
Let $H$ be a normal subgroup of a finite group $G$. Let ${\mathbb K}$ be a field with ${\rm char}\, {\mathbb K}$ not dividing $|G|$ and let $L$ be an irreducible ${\mathbb K}H$-module.  If for  all $x \not\in H$, the ${\mathbb K}H$-modules, $L_{H}$ and $x \otimes L$ are disjoint then the induced module $L^G$ of $G$ is irreducible.
\end{coro}

\noindent{\bf Proof}{   }Since $H^{(x)} = H$ in this case the proof is a direct consequence of Theorem 3.1. \hspace{.75in} $\Box$

\vspace{.05in}
\noindent
The above result is a particular case of the more general result given below. 
\vspace{.05in}

{\em Let $N$ be a normal subgroup of an arbitrary group G, let $R$ be a commutative ring and let $V$ be an irreducible $RN$-module. If $g \otimes V \not \cong V$ for all $g \in G\setminus N$, then $V^G$ is irreducible.} (See \cite{gk90}, p 96 for a proof.)
\vspace{.05in}

\begin{coro}
Let $H$ be a subgroup of a finite group $G$. Let ${\mathbb K}$ be a field with ${\rm char}\, {\mathbb K}$ not dividing $|G|$ and let $\rho$ be a one-dimensional representation of $H$. If for each $x \not \in H$, there exists  $y \in xHx^{-1}\cap H$ such that $\rho(y) \not = \rho(x^{-1}yx) $ then the induced monomial representation ${\rho}^G$ of $G$ is irreducible.
\end{coro}

\noindent{\bf Proof}{   }Let $L$ be a ${\mathbb K}H$-module which affords $\rho$. For any $x \in G$, let ${\rho}^x$ be the representation of $xHx^{-1}$ that is given by ${\rho}^x(g) = \rho(x^{-1}gx)$ for all $g \in xHx^{-1}$. Let us denote by $L^x$ the ${\mathbb K}xHx^{-1}$-module which affords ${\rho}^x$. Note that the underlying vector space for the module $L^x$ is $L$ itself. It is obvious that $L$ and $L^x$ are irreducible, as their dimension over ${\mathbb K}$ is 1.

It is given that for $x \not \in H$, there exists $y \in xHx^{-1}\cap H =: H^{(x)}$, such that $\rho(y) \not = \rho(x^{-1}yx)$. So there exists $y \in H^{(x)}$ such that 
$\rho(y) \not = \rho^{x}(y)$ or equivalently we have that for all $x \not \in H$, the one-dimensional representations $\rho$ and $\rho^x$ restricted to $H^{(x)}$ are not equal. Since they are one-dimensional representations, we get that they are not equivalent. In terms of modules, all this is saying is that the ${\mathbb K}H^{(x)}$-modules,  $L_{ H^{(x)}}$ and ${(L^x)}_{H^{(x)}}$ are non-isomorphic. It is easy to see that $x \otimes L$ and $L^x$ are isomorphic as ${\mathbb K}xHx^{-1}$-modules. Thus what we have is that the ${\mathbb K}H^{(x)}$-modules, $L_{ H^{(x)}}$ and ${(x \otimes L)}_{H^{(x)}}$ are not isomorphic. Since the modules involved are irreducible, being non-isomorphic is the same as being disjoint. So we have shown the following. 

For all $x\not \in H$, the ${\mathbb K}H^{(x)}$-modules, $L_{ H^{(x)}}$ and ${(x \otimes L)}_{H^{(x)}}$ are disjoint. So by Theorem 3.1, the induced module $L^G$ is irreducible as a ${\mathbb K}G$-module or equivalently the induced monomial representation ${\rho}^G$ of $G$ is irreducible. \hspace{.75 in} $\Box$

We end this section by mentioning a theorem which is a test for isomorphism of induced modules. (For a proof of this theorem see \cite{gk90} p 94.) This result will be used in the proof of the main theorem   presented in the next section.
\begin{theo}
Let $H_1$ and $H_2$ be subgroups of a finite group $G$ and let $L_i$ be ${\mathbb K}H_i$-modules for $i=1, 2$, where ${\mathbb K}$ is a field such that ${\rm char}\,{\mathbb K}$ does not divide $|G|$. Further let ${(L_i)}^G$ be irreducible as ${\mathbb K}G$-modules and let $H^{(x)} := x H_1x^{-1} \cap H_2$. Then ${(L_1)}^G$ and ${(L_2)}^G$ are not ${\mathbb K}G$-isomorphic if and only if, for all $x \in G$, the ${\mathbb K}H^{(x)}$-modules, $x \otimes L_1$ and $L_2$ are disjoint.
\end{theo}

\section{An adaptation of the Wigner--Mackey\\
method of little groups}

For any group $G$ and a field ${\mathbb K}$, let ${\rm Irr}_{\mathbb K}(G)$ denote the set of all irreducible representations of $G$ over ${\mathbb K}$ up to isomorphism. The set of all one dimensional representations of $ G$ over ${\mathbb K}$ forms an abelian group and will be denoted as ${\tilde {G}}$. 

Let ${N}$ be a normal subgroup of a group ${G}$. Then ${G}$ acts on ${\tilde {N}}$ as follows: given $\chi \in {\tilde {N}}$, $g \in {G}$, for all $a \in {N}$, we have ${\chi}^g(a) = \chi(g^{-1}ag)$. 

In this section we shall present a classification of irreducible representations of a finite group $G$ over a field ${\mathbb K}$ with respect to the conditions given below.
\begin{enumerate}
\item[(i)]${\rm char}\,{\mathbb K}$ does not divide $|G|$.
\item[(ii)]$G$ is a semidirect product of a normal abelian subgroup $N$ by a subgroup $H$.
\item[(iii)]All irreducible representations of $N$ over ${\mathbb K}$ have degree 1.
\end{enumerate} 

We also present a result which deals with constructing irreducible representations of $G$ but without the imposition of condition (iii). This result gives a classification of all irreducible ${\mathbb K}G$-modules which have a one-dimensional composition factor when restricted to $N$.

Let us assume that condition (ii) mentioned above holds.
We know that $G$ acts on ${\tilde N}$ as follows: given $\chi \in {\tilde N}$, $g \in G$, for all $a \in N$, we have ${\chi}^g(a) = \chi(g^{-1}ag)$. We shall denote by $I_{\chi}$ the stabiliser of ${\chi}$ in $G$. Note that since $N$ is abelian, we have that $N \leq I_{\chi}$. 

Let $H_{\chi} := I_{\chi} \cap H$. Then it is fairly obvious that $I_{\chi}$ is a semidirect product of $N$ by $H_{\chi}$. It can be shown quite easily that any $\chi \in {\tilde N}$ can be extended to a homomorphism from $I_{\chi}$ to ${\mathbb K}^*$, the multiplicative group of the field ${\mathbb K}$, in such a way that this extended homomorphism is the trivial map when restricted to $H_{\chi}$. So now we can regard $\chi$ as an element of $\tilde{I_{\chi}}$.

Further if $\rho$ is a representation of $H_{\chi}$ over ${\mathbb K}$ and the canonical projection of $I_{\chi}$ on $H_{\chi}$ is composed with $\rho$, then we get a representation of $I_{\chi}$. Thus we can form the tensor product representation $\chi \otimes \rho$ of $I_{\chi}$. It is easy to show that if $\rho$ is irreducible then the tensor product representation $\chi \otimes \rho$ is also irreducible and $\deg (\chi \otimes \rho) = \deg \rho$. 

We are now in a position to state and prove our classification theorems.

\begin{theo}
Let $G$ be a finite group which is a semi-direct product of an abelian group $N$ by a subgroup $H$. Let ${\mathbb K}$ be a field such that ${\rm char} \, {\mathbb K}$ does not divide $|G|$ and let all irreducible representations of $N$ over ${\mathbb K}$ be of degree 1. Let $O_1, \ldots, O_t$ be the distinct orbits under the action of $G$ on ${\tilde N}$ and let ${\chi}_j$ be a representative of the orbit $O_j$. Let $I_j$ denote the stabiliser of ${\chi}_j$ and let $H_j := I_j \cap H$. For any irreducible representation $\rho$ of $H_j$, let ${\theta}_{j, \rho}$ denote the representation of $G$ induced from the irreducible representation ${\chi}_j \otimes \rho$ of $I_j$. Then
\begin{enumerate}
\item[\rm{(i)}]${\theta}_{j, \rho}$ is irreducible.
\item[\rm{(ii)}]If ${\theta}_{j, \rho}$ and ${\theta}_{j', {\rho}'}$ are isomorphic, then $j = j'$ and $\rho$ is isomorphic to ${\rho}'$.
\item[\rm{(iii)}]Every irreducible representation of $G$ is isomorphic to one of the ${\theta}_{j, \rho}$.
\end{enumerate}
\end{theo}

\noindent{\bf Proof}{   }Let $W_{\rho}$ be a ${\mathbb K}H_j$-module affording the representation $\rho$. Let $W_{j,\rho}:= V_j \otimes W_{\rho}$ where $V_j$ is a one-dimensional representation space for $I_j$ affording the character $\chi_j$. Then we may regard $W_{j,\rho}$ as a ${\mathbb K}I_j$-module affording the irreducible tensor product representation $\chi_j \otimes \rho$. So ${(W_{j, \rho})}^G$ is a ${\mathbb K}G$-module affording ${\theta}_{j, \rho}$. 

For part (i) it is sufficient to show that ${(W_{j, \rho})}^G$ is irreducible as a ${\mathbb K}G$-module. By Theorem 3.1, it is sufficient to show that for any $x \not \in I_j$ and for ${I_j}^{(x)}:= xI_jx^{-1} \cap I_j$, the ${\mathbb K}{I_j}^{(x)}$-modules $W_{j,\rho}$ and $x \otimes W_{j,\rho}$ are disjoint.

Suppose for some $x \not \in I_j$, we have that the ${\mathbb K}{I_j}^{(x)}$-modules $W_{j,\rho}$ and $x \otimes W_{j,\rho}$, have a composition factor in common. Then if we further restrict these modules to the subgroup $N$ of ${I_j}^{(x)}$, we must have that as ${\mathbb K}N$-modules, $W_{j,\rho}$ and $x \otimes W_{j,\rho}$ have a composition factor in common.

Now for any $w \in W_{j,\rho}$ and $a \in N$, let $a \cdot w$ denote the action of $a$ on $w$. Then it is easy to see that
$a \cdot w =\chi_j(a) \, w $. So it is obvious that irreducible ${\mathbb K}N$-submodules of the ${\mathbb K}N$-module $W_{j,\rho}$ are one dimensional and afford the character $\chi_j$. Similarly we can show that irreducible ${\mathbb K}N$-submodules of the ${\mathbb K}N$-module $x \otimes W_{j,\rho}$ are one dimensional and afford the character ${\chi_j}^x$.

Since $x \not \in I_j$, we have that $\chi_j \not = {\chi_j}^x$. So the ${\mathbb K}N$-modules, $W_{j,\rho}$ and $x \otimes W_{j,\rho}$ cannot have a composition factor in common. Therefore our assumption that the ${\mathbb K}{I_j}^{(x)}$-modules $W_{j,\rho}$ and $x \otimes W_{j,\rho}$, have a composition factor in common is false. Hence by Theorem 3.1 we have that ${(W_{j,\rho})}^G$ is irreducible or equivalently that $\theta_{j, \rho}$ is irreducible.

For part (ii), let us assume that ${\theta}_{j, \rho}$ and ${\theta}_{j', {\rho}'}$ are isomorphic. As in the previous part, we shall assume that the ${\mathbb K}I_j$-module $W_{j,\rho}$ and the ${\mathbb K}I_{j'}$-module $W_{j',{\rho}'}$ afford the representations ${\chi}_j \otimes \rho$ and ${\chi}_{j'} \otimes {\rho}'$ respectively. We are given that the irreducible induced modules ${(W_{j,\rho})}^G$ and ${(W_{j',{\rho}'})}^G$ are isomorphic as ${\mathbb K}G$-modules. 

By Theorem 3.2, we have that for some $x \in G$, the ${\mathbb K}(x I_j x^{-1} \cap I_{j'})$-modules $x \otimes W_{j,\rho}$ and $W_{j',{\rho}'}$ have a composition factor in common. By restriction we can consider both these modules as ${\mathbb K}N$-modules and they will still have a composition factor in common. But then as in the previous part, we shall get that ${{\chi}_j}^x = {\chi}_{j'}$. Thus ${\chi}_j$ and ${\chi}_{j'}$ are in the same orbit and so we have $j = j'$.

Since $j=j'$, we have that ${{\chi}_j}^x = {\chi}_{j}$. So $x \in I_j$ and $x I_j x^{-1} \cap I_{j'} = I_j$. Hence we now have that the ${\mathbb K}I_j$-modules $W_{j,\rho}$ and $W_{j,{\rho}'}$ have a composition factor in common. But these modules are also irreducible as ${\mathbb K}I_j$-modules. Thus we get that $W_{j,\rho}$ and $W_{j,{\rho}'}$ are isomorphic as ${\mathbb K}I_j$-modules. Since $H_j \leq I_j$, if we restrict $W_{j,\rho}$ and $W_{j,{\rho}'}$ to $H_j$, then they will still be isomorphic as ${\mathbb K}H_j$-modules. But $W_{j,\rho}$ restricted to $H_j$ affords the representation $\rho$ and $W_{j,{\rho}'}$ restricted to $H_j$ affords the representation ${\rho}'$. Thus $\rho \cong {\rho}'$.

For the last part, let us assume that $V$ is an irreducible ${\mathbb K}G$-module with representation $\theta$. Now $V_N$ is completely reducible and since all irreducible representations of $N$ over ${\mathbb K}$ are of degree 1, the abelian group, ${\tilde N}$, is the set of all irreducible representations of $N$ over ${\mathbb K}$. So ${\tilde N} = {\rm Irr}_{\mathbb K}(N)$ and we can write
$$
V_N = \oplus_{\chi \in {\tilde N}}V_{\chi}
$$ 
where $V_{\chi} = \{v \in V \mid \theta(a)v = \chi(a) v, \ \forall a \in N \}$.

For any $x \in G$ and $\chi \in {\tilde N}$, we have $\theta(x) (V_{\chi}) = V_{\chi^x}$. Since $V \not = 0$, there exists $\chi$, such that $V_{\chi} \not = 0$. We can assume without loss of generality, there exists $j$, such that $\chi = \chi_j$. 

Now for any $x \in I_j$, we have $\theta(x) (V_{\chi_j}) = V_{{\chi_j}^x} = V_{\chi_j}$. So $V_j := V_{\chi_j}$ is a ${\mathbb K}I_j$-submodule of $V$ and is completely reducible as a ${\mathbb K}I_j$-module. Let $W$ be any irreducible ${\mathbb K}I_j$-submodule of $V_j$. So for any $a \in N$ and $w \in W$, we have $\theta(a) w = \chi_j(a) w$. By restriction we can regard $W$ as a ${\mathbb K}H_j$-module. It is easy to show that any ${\mathbb K}H_j$-submodule of $W_{H_j}$ is also a ${\mathbb K}I_j$-submodule of $W$. Thus $W_{H_j}$ is irreducible. 

Let $\rho$ be the irreducible representation afforded by $W_{H_j}$. Then for any $w \in W$ and $h \in H_j$, we have $\rho(h) w = \theta(h) w$. Further, if $V_j$ is a one-dimensional representation space for $I_j$ affording the character $\chi_j$ then, as in the previous parts, we shall assume that $W_{j, \rho} := V_j \otimes W$ is a ${\mathbb K}I_j$-module affording the representation $\chi_j \otimes \rho$.

Now let $g = ah \in I_j$ where $a \in N$ and $h \in H_j$. Then for $w \in W$ and $v \in V_j$, 
\begin{eqnarray*}
({\chi}_j \otimes \rho \, (g)) (v \otimes w) & = & \chi_j(a) v \otimes  \rho (h) (w) \\
					   & = & \chi_j(a) v \otimes(\theta(h)(w))\\
                                 & = &  v \otimes \chi_j(a)(\theta(h)(w))\\
					   & = & v \otimes\theta(a)(\theta(h) (w))\\
					   & = & v \otimes\theta(ah) (w )\\
					   & = & v \otimes \theta(g)(w) \: .
\end{eqnarray*} 
Thus $W \cong W_{j, \rho}$ as ${\mathbb K}I_j$-modules and we have shown  that $V_{I_j}$ has a composition factor isomorphic to $W_{j, \rho}$, which affords the representation $\chi_j \otimes \rho$. So by Remark 2.1, we have $i(V_{I_j}, W_{j, \rho}) \not = 0$ and by FRT we get that $i(V, {(W_{j, \rho})}^G )\not = 0$. But $V$ is an irreducible ${\mathbb K}G$-module and by part (i), so is ${(W_{j, \rho})}^G$. Since ${\rm Hom}_{{\mathbb K}G}(V, {(W_{j, \rho})}^G) \not = 0$, by Schur's Lemma we get that $V$ is isomorphic to ${(W_{j, \rho})}^G$ as ${\mathbb K}G$-modules. Equivalently we have that $\theta \cong \theta_{j, \rho}$. \hspace{.75in} $\Box$

Our next result is similar to the above result except that we no longer impose the condition that all irreducible ${\mathbb K}N$-modules have dimension 1. 

\begin{theo}
Let $G$ be a finite group which is a semi-direct product of an abelian group $N$ by a subgroup $H$. Let ${\mathbb K}$ be a field such that ${\rm char} \, {\mathbb K}$ does not divide $|G|$. Let $O_1, \ldots, O_t$ be the distinct orbits under the action of $G$ on ${\tilde N}$ and let ${\chi}_j$ be a representative of the orbit $O_j$. Let $I_j$ denote the stabiliser of ${\chi}_j$ and let $H_j := I_j \cap H$. For any irreducible representation $\rho$ of $H_j$, let ${\theta}_{j, \rho}$ denote the representation of $G$ induced from the irreducible representation ${\chi}_j \otimes \rho$ of $I_j$. Then
\begin{enumerate}
\item[\rm{(i)}]${\theta}_{j, \rho}$ is irreducible.
\item[\rm{(ii)}]If ${\theta}_{j, \rho}$ and ${\theta}_{j', {\rho}'}$ are isomorphic, then $j = j'$ and $\rho$ is isomorphic to ${\rho}'$.
\item[\rm{(iii)}]Let $V$ be an irreducible ${\mathbb K}G$-module affording the representation $\theta$ and let $V_N$ have a one-dimensional composition factor. Then $\theta$ is isomorphic to one of the ${\theta}_{j, \rho}$.
\end{enumerate}
\end{theo}

\noindent{\bf Proof}{   }The proofs of parts (i) and (ii) are identical to the proofs of the same parts in Theorem 4.1.

For part (iii), let us assume that $V$ is an irreducible ${\mathbb K}G$-module affording the representation $\theta$ and let $V_N$ have a one-dimensional composition factor. So there exists a ${\mathbb K}N$-submodule, $V_o$ of $V$ and $\chi_o \in {\tilde N}$ such that $\theta = \chi_o$ on $V_o$.

For any $\chi \in {\tilde N}$, let $V_{\chi} = \{v \in V \mid \theta(a)v = \chi(a) v, \ \forall a \in N \}$. For any $x \in G$ and $\chi \in {\tilde N}$, we have $\theta(x) (V_{\chi}) = V_{\chi^x}$.

It is clear that $V_o \subseteq V_{\chi_o}$ and so $V_{\chi_o} \not = 0$. Further by the above paragraph we may assume that $\chi_o = \chi_j$ for some $j$. 

Thus we have that the subspace $V_{\chi_j}$ of $V$ is non-zero. As in the proof of part (iii) of Theorem 4.1, we can regard $V_{\chi_j}$ as a ${\mathbb K}I_j$-module. Proceeding exactly as we did in the proof from that point, we can show that $V_{I_j}$ has a composition factor affording the representation $\chi_j \otimes \rho$ for some irreducible representation $\rho$ of $H_j$ and that $\theta$ is isomorphic to $\theta_{j, \rho}$. \hspace{.75in} $\Box$

We conclude with a note on some examples of groups and fields that satisfy the hypothesis of Theorem 4.1.
\begin{note}
Let $G$ be a finite group which is a semidirect product of a normal abelian subgroup $N$ of exponent $m$ by a subgroup $H$. Let ${\mathbb K}$ be a finite field with $q$ elements such that ${\rm char}\, {\mathbb K}$ does not divide $|G|$ and such that $m$ divides $q-1$. Then all irreducible ${\mathbb K}N$-modules have dimension 1 and so $G$ satisfies the hypothesis of Theorem 4.1.
\end{note}

\bibliography{vision1}
\end{document}